\documentclass[12pt]{article}
\usepackage{amssymb}
\usepackage{amsmath}
\usepackage{amsthm}
\usepackage{graphicx}
\usepackage{graphics}
\usepackage{amsthm}

\input pictex.tex
\theoremstyle{definition}

\newtheorem{thm}{Theorem}[section]

\newtheorem{prop}[thm]{Proposition}

\newtheorem{lem}[thm]{Lemma}

\newtheorem{rem}[thm]{Remark}

\newtheorem{defn}[thm]{Definition}

\newtheorem{ex}[thm]{Example}

\newtheorem{obs}[thm]{Observation}

\newcommand{\Diff}{\mathrm{Diff}}
\newcommand{\Homeo}{\mathrm{Homeo}}
\newcommand{\PL}{\mathrm{PL}}

\input epsf

\begin{document}
\date{}
\author{Leonardo Dinamarca Opazo}

\title{Distortion in groups of generalized piecewise-linear transformations}
\maketitle

\vspace{-0.4cm}

\noindent{\bf Abstract.} For each natural number $n$, we consider the subgroup $\mathcal{R}_n$ of $\Homeo_+[0,1]$ made by the elements 
that are linear except for a subset whose Cantor-Bendixson rank is less than or equal to $n$. These groups of generalized piecewise-linear 
transformations yield an ascending chain of groups as we increase $n$. We study how the notion of distorted element changes along this chain. 
Our main result establishes that for each natural number $n$, there exits an element that is undistorted of $\mathcal{R}_n$ yet distorted in 
$\mathcal{R}_{n+1}$. Actually, such an element is explicitly constructed.

\vspace{0.2cm}


\vspace{0.2cm}


%%%%%%%%%%%%%%%%%%%%%%%%%%%%%%%%%%%%%%%%%%%%%%%%%%%%%%%%%%%%%%%%%%%%%%

\section{Introduction}

The notion of distorted element was introduced by Mikhail Gromov \cite{gromov}. Given a finitely generated group $\Gamma$, we fix a finite system 
of generators $S$, and we denote $l_S(\cdot)$ the corresponding word-length. An element $f \in \Gamma$ is said to be distorted if 
$$\lim_{n \to \infty} \frac{l_S(f^n)}{n} = 0.$$
(Note that this condition does not depend on the choice of the finite generating system.) Given an arbitrary group $G$, an element 
$f \in G$ is said to be distorted if there exists a finitely generated subgroup $\Gamma \subset G$ containing $f$ so that 
$f$ is distorted in $\Gamma$ in the sense above. 

This notion has been studied in different groups, for example in the automorphisms group of a shift  \cite{CS, CFKP}, in 
the Cremona group \cite{Can} and mainly in groups of transformations \cite{GL, hur, Le roux Kat, polt}. We will focus in the 
last direction. In \cite{FH}, John Franks and Michael Handel asked whether irrational rotations are distorted in $\Diff^1_+(\mathbb{S}^1)$.  
Danny Calegari and Michael Freedman (see \cite{CF}) gave a positive answer to this.

\begin{thm} 
    {\em Irrational rotations are distorted in $\Diff^1 _+(\mathbb{S}^1)$.}
\end{thm}

Later, Arthur Avila extended this result to infinite differentiability (see \cite{Avila}). To state his result in full generality, 
recall that $f$ is said \textit{recurrent} if $\liminf_{n\rightarrow\infty} d(f^n,id) =0,$ where $d$ is a metric inducing the 
$C^{\infty}$ topology on $\Diff_+^{\infty}(\mathbb{S}^1).$

\begin{thm} \label{Av}
    {\em  In $ \Diff_+^{\infty}(\mathbb{S}^1)$, every recurrent element is distorted.}
\end{thm}

 Emmanuel Militon generalized  this result to any compact manifold (see \cite{Mil}). The following related question remains open.

 \vspace{0.2 cm}
    {\em  \noindent{\textbf{Question 1.}} Are irrational rotations on the circle distorted in the group of real-analytic difeomorphisms $\Diff_+^{\omega}(\mathbb{S}^1)$?}
 \vspace{0.2 cm}

 In \cite{Na1}, Andrés Navas gave a list of concrete questions about group actions on $1$-dimensional manifolds. One of these questions was  
 the existence of distorted elements in $\PL_+(\mathbb{S}^1)$, the group of piecewise-linear homeomorphisms of the circle. This question was 
 answered affirmatively by Juliusz Banecki and Tomasz Szarek in \cite{BS}; also, their proof can be adapted to obtain an alternative proof of 
 Avila's result  (Theorem \ref{Av}). 
 
\begin{thm}
    {\em  All irrational rotations are distorted in $\PL_+(\mathbb{S}^1)$.}
\end{thm}

 There exist many interesting countable subgroups of $\PL_+(\mathbb{S}^1)$ for which it is unkknow whether they contain distorted elements 
 or not. For instance, let $F_\tau$ be the Thompson group (introduced by Sean Cleary in \cite{Clea}) associated to the golden number $\tau$, defined 
 as the set of elements in $\PL_+(\mathbb{S}^1)$ with break points in $\mathbb{Z}[\tau]$ and slopes that are powers of $\tau$. The next question 
 was suggested by Yash Lodha.

\vspace{0.2cm}
    {\em \noindent{\textbf{Question 2.}} Does $F_\tau$ contain distorted elements?}
\vspace{0.2cm}

The situation in $\PL_+(I)$, the group of piecewise-linear homeomorphisms of the interval $I=[0,1]$ is different, since the unique distorted element is the identity. 
To see this, let $L:\PL_+(I)\rightarrow \mathbb{R}$ be the function defined as 
$$L(f):=\log \max_{x\in I}\{D_Rf (x),D_Rf^{-1}(x)\}\mbox{ for all } f\in \PL_+(I).$$
 Let us note that for $f,g \in \PL_+(\mathbb{S}^1)$ the chain rule implies that 
 $$L(fg)\leq L(f)+L(g).$$
 If $f$ is a nontrivial element, then we have that $ L (f^n) \geq C\cdot n$ for some positive constant $C.$ 
 This easily implies that $f$ is undistorted (see the Criterion in Section \ref{section exam}). 
 
 Since there is no nontrivial distorted element in $\PL_+(I)$ we will relax the condition on the quantity of break points to consider larger groups.

Let $n$ be a natural number. Following a suggestion of Lodha, we  define the \textit{\textbf{group of generalized piecewise-linear transformations}},  
that we will denote by $\mathcal{R}_n$ (see Definition \ref{gp Rn}), as the subgroup of $\Homeo_+[0,1]$ made by the elements that are linear 
except for a subset of $[0,1]$ (that we call set of break points) whose Cantor-Bendixson rank is smaller than or equal to $n$. The fact that 
this is actually a group is not totally obvious (see Proposition \ref{gp}). Once this established, we get an ascending chain of groups: 
$$  \PL_+(I)=\mathcal{R}_0\leq\mathcal{R}_1\leq\mathcal{R}_2\leq\mathcal{R}_3\leq\dots$$ 

A general problem is to study how the notion of distortion change along a filtration of groups. In \cite{Na}, He addressed the following question:

\vspace{0.2 cm}

{\em \noindent{\textbf{Question 3.}} \label{qs andres} Given $1 \leq r < s$, does there exist an undistorted element $f \in \Diff^{s}_+ (M)$ that is distorted when considered 
as an element of $\Diff^{r}_+ (M)$~?}   

\vspace{0.2cm}

He proved that this is the case for $M=[0,1]$, $r = 1$ and $s=2$. Actually, in his example, undistortion holds in the larger group 
$\Diff^{1+bv}_+([0,1])$ of $C^1$ diffeomorphisms with derivative of bounded variation. Moreover, in \cite{DE}, finer computations 
concerning intermediate regularity are made and. As an outcome, it is established that that there exists a diffeomorphism that is 
distorted in $\Diff_+^{1+\alpha}[0,1]$ but undistorted in $\Diff_+^{2}[0,1]$. However, for the case of the closed interval, the question 
above remains open in regularity larger than $2$ and it is completely open for higher dimensional manifolds. 

In the same spirit, Lodha asked the following question: given $n$ a natural number, does there exist an undistorted element in 
$\mathcal{R}_n$ that is distorted in  $\mathcal{R}_{n+1}$? The main result of this work answers this question positively.

\vspace{0.5cm}

\noindent{\textbf{Theorem A}}
{\em For each natural number $n$, there exists $f\in Homeo_+[0,1]$ that is undistorted in $ \mathcal{R}_n$} but distorted in $\mathcal{R}_{n+1}$.

\vspace{0.5cm}

The strategy to prove this theorem is as follows: given a natural number $n$, we define a length function in $\mathcal{R}_n$ that we will use 
to detect undistorted elements. Then, we build an element that we force to be undistorted (using the length function) and then, we establish  
distortion in $\mathcal{R}_{n+1}$ via a combination of the diagonal trick and the Mather's argument. 

 The groups $\mathcal{R}_n$ remain mysterious to us. 
 A study from the perspective of large scale geometry in the sense of Rosendal \cite{ros} would be a good starting point.

\section{Preliminaries}
In this section, we provide the definition of distorted elements in groups, discuss some of their properties, and present several basic examples.

\begin{defn}
{\em Let $\Gamma =\langle S\rangle$ be a finitely generated group. Given an element $f\in \Gamma$, we define the \textbf{word-length} of $f$ with respect to $S$ as the minimum integer $n$ such that $f$ can be written as product of $n$ elements of $S.$ We use the following notation:}   
$$l_S(f)=||f||:=\min\{n\,|\,\,f=s_1\cdot s_2\cdot ...\cdot s_n,\,\mbox{ with } s_1,...,s_n \in S \}.$$ 
     
\end{defn}

 An important property satisfied by the word-length function $l_S:\Gamma \rightarrow \mathbb{N}$
 is the subadditive inequality: for every pair of elements $f,g \in \Gamma$, we have
 $$l_S(gh)\leq l_S(g)+l_S(h).$$

   \begin{rem}
        The word metric is given by $d_S(f,g):=l_S(f^{-1}g)$. The topology induced by the word metric is the discrete topology for any generating set.
   \end{rem}

Observe that if we take two generating sets of the group, we obtain a lipschitz relation between the induced word-length functions. Indeed, if $\Gamma=\langle S\rangle=\langle T\rangle$, then for each $g\in\Gamma$, the following inequality holds
\begin{equation}\label{lip rel} 
    \frac{1}{C'}\cdot l_S(g)\leq l_T(g)\leq C \cdot l_S(g),
\end{equation}
where $C=\max_{s\in S}l_T(s) $ and $C'=\max_{t\in T}l_S(t)$.

     We next define an important limit for the rest of the text. Let $f\in \Gamma=\langle S\rangle $ be any element. By the subadditive inequality, we have
       \begin{equation*}
        l_S(f^{n+m})\leq l_S(f^n)+l_S(f^m),  
       \end{equation*}
     for all natural numbers $n$ and $m$. This means that the sequence $(l_S(f^n))$ is subadditive, therefore by Fekete's Lemma \cite{Fek}, the following limit exists: 
        $$\lim_{n\rightarrow \infty} \frac{l_S(f^n)}{n}=\inf_{n\geq 1}\frac{l_S(f^n)}{n}.$$
\begin{rem}
        Note that if the limit is equal to $0$, then the value does not depend on the choice of the finite generating system (because of (\ref{lip rel})).
        \end{rem}
       
       Now we define the concept of distortion in groups, which is the crucial notion of this work. 
        \begin{defn}
            {\em Let $\Gamma$ be a finitely generated group. We say that $f\in \Gamma$ is  \textbf{distorted} if the limit above is equal to $0.$ 
            If the limit is different from $0$ we say that $f$ is \textit{undistorted}. Note that being a distorted element is a property of the group and the element.}
        \end{defn}

           Given a non-torsion distorted element $f$ in a group $G=\langle S\rangle$, we will pay attention to the growth rate of $l_S(f^n)$ as a function of $n$ to analyze 
           different behavior. The \textbf{distortion function} $D_{S,f}:\mathbb{N}\rightarrow\mathbb{N}$ is defined by  
\begin{align*}
    D_{S,f}(n):=\max\{k : l_S(f^k)\leq n\}.
\end{align*}

 We have seen that the notion of distorted element does not depend of $S$. 
 In order to have the same property for the distortion function, we consider an equivalence relation. Given two nondecresing functions
 \begin{align*}
     h,g:\mathbb{N}\rightarrow\mathbb{N},
 \end{align*}
 we write $h\preceq g$ if there exists $k\geq 1$ such that
 \begin{align*}
     h(n)\leq k\cdot g(k\cdot n+k)+k \mbox{ for all } n\in \mathbb{N}
 \end{align*}
 and we write $h\sim g$ if $h\preceq g$ and $g\preceq h$. We have that $\sim$ is an equivalence relation. 
 With this definition, if $S'$ is another generating system of $G$, we obtain  that $D_{S,f}\sim D_{S',f}$ and 
 then we simply denote the distortion function of $f$ as $D_f.$

\begin{defn}
        {\em For a general group $G$, we say that an element is distorted if there exists a finitely generated subgroup of $G$ such that the element is distorted 
        inside this finitely generated subgroup. }
\end{defn} 

\subsection{Examples of distorted and undistorted elements} \label{section exam}

In this section we will give examples of distorted and undistorted elements in groups. First, we recall the notion of a length function, and we define 
those that we will use to state the criterion that we apply to show that an element is undistorted. Such length functions arise from geometric group theory, 
dynamical systems, algebraic geometry and many other branches of mathematics (see \cite{ye} for a survey). After this, we present several 
examples of distorted and undistorted elements 

\begin{defn}\label{defn length}
   {\em Let $G$ be a group. A \textbf{length function} on $G$ is a function $L:G \rightarrow \mathbb{R}\cup \{0\}$ such that $L(Id) = 0$, $L(g)=L(g^{-1})$ and $L(gh)\leq L(g)+L(h)$ for all $g,h\in G$.}
   \end{defn}

Given a length function $L$ on a group $G$ and an element $f\in G$, we define the \textit{\textbf{stabilazed limit}} as 
$$\lim_{n\rightarrow\infty}\frac{L(f^n)}{n}.$$
Again, existence of this limit follows from Fekete's lemma. 

\begin{lem} \label{lenght}
   {\em Let $G$ be a group and $L$ a length function on $G$. If $g\in G$ is distorted, then the stabilazed limit is equal zero}
\end{lem}
\noindent{{\bf Proof.}} Let $g$ be distorted in $ \Gamma=\langle S\rangle$, a finitely generated generated subgroup of $G$. 
We let $C:=\max\{L(s):\,s\in S\}$. Then, by subadditivity of $L$, we have
\vspace{0.5cm}
\begin{equation*}
    \lim_{n\rightarrow\infty}\frac{L(f^n)}{n}\leq C\cdot\lim_{n\rightarrow\infty}\frac{l_S(f^n)}{n}=0.      
\end{equation*}

\vspace{0.3cm}

We use the previous lemma to test that an element is undistorted. 
%Our Criterion is the lemma restated. It say that if we build a length function in a group $G$ such that the stabilized limit is different from zero then the element cannot be distorted.

\vspace{0.5cm}
\textbf{Criterion:} \textit{If there exists a length function $L:G\rightarrow\mathbb{R}$ such that for an element $f\in G$ the stabilazed limit is positive, then $f$ is undistorted.} 

\vspace{0.5cm}
The first example of distortion concerns the \textbf{Baumslag-Solitar} group and its realization as a group of homeomorphisms of the real line.
\begin{ex} \label{bs}
         Let $f$ and $g$ be the elements in $\mbox{Homeo}_+(\mathbb{R})$ defined by 
        $f(x)=x+1$ and $g(x)=2x$. Consider the group $\Gamma=\langle S\rangle$ with $ S:=\{ f^{\pm 1},g^{\pm 1}\}$. An explicit computation shows that 
        $$g^nfg^{-n}(x)=f^{2^n}(x) \mbox{ for all } n\geq 1.$$
     Hence, $f$ is distorted in $\langle S\rangle$, since
        $$\lim_{n\rightarrow\infty} \frac{l_S(f^{2^n})}{2^{n}}\leq\lim_{n\rightarrow\infty}\frac{2n+1}{2^{n}}=0.$$
    
\end{ex}

In the next example, we exhibit distorted and undistorted element in a finite chain of groups. We will see explicitly how an element is undistorted in a group but distorted in a larger one. We consider a group of upper triangle matrices of $5\times5$ with integer coefficients and a chain of three subgroups.   

\begin{ex} \label{exam dis}
    Let us consider the group
\begin{align*}
            \mathcal{H}_5:=\left\{\begin{pmatrix}
    1 & a_{1,2} & a_{1,3} &a_{1,4} & a_{1,5}\\
    0 & 1 & a_{2,3} & a_{2,4} & a_{2,5}\\
    0 &  0 & 1 & a_{3,4} &a_{3,5}\\
    0 & 0 &  0 & 1 & a_{4,5}\\
    0 & 0 & 0 & 0 &1
    \end{pmatrix}
    : a_{i,j}\in \mathbb{Z}\right\}.
     \end{align*}
Also, we consider the subgroups of $\mathcal{H}_5$ defined by
\begin{align*}
    \Gamma_1 &:=\{A\in \mathcal{H}_5: a_{1,2}=a_{1,3}=a_{1,4}=a_{1,5}=a_{2,3}=a_{2,4}=0\},\\
     \Gamma_2 &:=\{A\in \mathcal{H}_5: a_{1,2}=a_{1,3}=a_{1,4}=0\}.
\end{align*}
We have a finite chain of groups
$$ \Gamma_1\leq \Gamma_2\leq \mathcal{H}_5.$$
Let us define the matrices in $\mathcal{H}_5$ given by the formula
$$E_{i,j}:=e_{i,j}+Id,$$
where $2\leq j\leq 5$, $1\leq i < j$ and $e_{i,j}$ are the canonical matrices. Observe that the set $S$ made by the matrices $E_{i,j}$ is a finite generating system of $\mathcal{H}_5$. 
Moreover, for $k\in \{1,2\}$, we have that $S_k=S\cap\Gamma_k$ is a finite generating system of $\Gamma_k$. An explicit computation shows that
$$E_{3,5}^{n^2}=[E_{3,4}^n,E_{4,5}^n],$$    
so
$$\lim_{n\rightarrow\infty}\frac{l_{S_1}(E_{3,5}^{n^2})}{n^2}\leq \lim_{n\rightarrow \infty} \frac{4n}{n^2}=0.$$
This means that $E_{3,5}$ is distorted in $\Gamma_1.$ On the other hand, the element $E_{2,5}$ is undistorted in $\Gamma_1$. 
To show this, it suffices to consider $L_1:\Gamma_1\rightarrow\mathbb{N}$ defined as $$L_1(A):=|(A)_{2,5}|$$
for each $A\in \Gamma_1$. Note that $L_1$ is a length function and $L_1(E_{2,5}^n)=n$. Then, our Criterion assert that $E_{2,5}$ is undistorted in $\Gamma_1$. However, $E_{2,5}$ is distorted in $\Gamma_2$. Actually, we have that 
$$E_{2,5}^{n^2}=[E_{2,3},E_{3,5}^{n^2}]=[E_{2,3},[E_{3,4}^n,E_{4,5}^n]],$$
so
$$\lim_{n\rightarrow\infty}\frac{l_{S_2}(E_{2,5}^{n^2})}{n^2}\leq \lim_{n\rightarrow \infty} \frac{8n+2}{n^2}=0.$$
On the other hand, we have that $E_{1,5}$ is undistorted in $\Gamma_2$. To check this, 
we just repeat the argument before with the length function $L_2:\Gamma_2\rightarrow\mathbb{N}$ 
defined as $L_2(A):=|(A)_{1,5}|$ for each $A\in \Gamma_2$. Nevertheless, we claim that 
$E_{1,5}$ is distorted in $\mathcal{H}_5$. Indeed, we have that 
$$E_{1,5}^{n^4}=[E_{1,3}^{n^2},E_{3,5}^{n^2}]=[[E_{1,2}^n,E_{2,3}^n],[E_{3,4}^n,E_{4,5}^n]],$$
and so
$$\lim_{n\rightarrow\infty}\frac{l_{S}(E_{1,5}^{n^4})}{n^4}\leq \lim_{n\rightarrow \infty} \frac{16n}{n^4}=0.$$
\end{ex}

\vspace{0.2cm}

 \begin{rem}
In Example \ref{exam dis}, we have that the elements $E_{3,5}$ and $E_{2,5}$ have distortion function of order $g(n)=n^2$ and the element $E_{1,5}$ 
has distortion function of order $g(n)=n^4$. One can increase the size of the matrix to produce a longer finite chain of subgroups with distorted elements 
with polynomial distortion function of degree as large as we want. In Example \ref{bs}, the element $f$ has distortion function of order $g(n)=2^n.$ 
\end{rem}

    Let $f:\mathbb{N}\rightarrow\mathbb{R}$ be a function. We say that $f$ is {\bf superpolynomial} if 
    \begin{align*}
        \lim_{n\rightarrow\infty} \frac{\log (f(n))}{\log(n)}=\infty.
    \end{align*}
    We say that $f$ is {\bf subexponential} if
    \begin{align*}
        \lim_{n\rightarrow\infty} \frac{\log (f(n))}{n}=0.
    \end{align*}
    Finally, $f$ is said to have {\bf intermediate growth} if $f$ is both subexponential and superpolinomial. We are not able to exhibit an exmaple of an element with distortion function of intermediate growth, eventhouth for specialists it must be known.

The next example can be considered as an affirmative answer to Question $3$ for the case of $M=\mathbb{S}^1$ and regularity $C^0$ versus $C^1.$

\begin{ex}
    In the group Homeo$_+(\mathbb{S}^1)$, Calegari and Freedman prove that every element is distorted 
    (see \cite{CF}). If we go to $C^1$ regularity, we observe that the existence of an hyperbolic fixed point is an obstruction to $C^1$-distortion. 

 Indeed, define $L:\Diff_+^1(\mathbb{S}^1)\rightarrow\mathbb{R}_+$ by
 $$L(f):=\log \max\{ ||Df||_{\infty},||Df^{-1}||_{\infty}\}.$$
Note that $L$ is a length function. Now, if $f\in \Diff_+^1(\mathbb{S}^1)$ has an hyperbolic fixed point $x$, 
i.e, $f(x)=x$ and  $\lambda := \log Df(x) \neq 0$, then, by the chain rule,
$$L(f^n)\geq n|\lambda|.$$
This implies that $f$ is undistorted in $\Diff_+^1(\mathbb{S}^1)$.
\end{ex}
 
 \subsection{Cantor-Bendixson rank}
 
In this section we recall the notion of Cantor-Bendixson rank and state the properties that will be useful for us. 
We will use this rank to define the groups $\mathcal{R}_n$ (see Definition \ref{gp Rn}). Given a subset $A$ of 
the interval $[0,1]$, we denote by $A'$ its \textit{derived set}, i.e, the set of accumulation points of $A$. 
 
\begin{defn}
\em{ Let $X\subset [0,1]$ be a set. The \textbf{Cantor-Bendixson rank}, denoted by $r(X)$, is defined by induction as follows:
\begin{itemize}
    \item If $X$ is a finite set, then $r(X)=0$.
    \item $r(X)=r(X')+1.$
\end{itemize}
In this work, we only consider subset of the interval with finite Cantor-Bendixson.}
\end{defn}

\begin{prop} \label{prop rank}
\em{Let $X,Y\subseteq[0,1]$ and let $g\in$ Homeo$_+[0,1]$. Then 
\begin{itemize}
    \item $r(X)=r(g(X))$.
    \item If $r(X)\leq n$ and $r(Y)\leq n$ for a certain $n \in \mathbb{N}$, then $r(X \cup g(Y))\leq n.$
\end{itemize}}
\end{prop}
\noindent{\bf Proof.} We prove both statements by induction. For the first statement, the base case $(r(X)=0)$ holds because a 
homeomorphism sends a finite set into a finite set. Suppose that $r(X)=n$. By the induction hypothesis, $r(g(X'))=r(X')=n-1,$ hence 
$$r(g(X))=r(g(X)')+1=r(g(X'))+1=(n-1)+1=n.$$

For the second statement, the case $n=0$ means that both $X$ and $Y$ are finite sets, hence $X\cup g(Y)$ is a finite set, therefore, $r(X \cup g(Y))=0$. Suppose $r(X)\leq n$ and $r(Y)\leq n$. Observe that
$$(X \cup g(Y))'=X' \cup g(Y)'=X' \cup g(Y').$$
Thus, by the induction hypothesis, it follows that $r((X \cup g(Y))')=r(X' \cup g(Y'))\leq n-1$ because $r(X')\leq n-1$ and $r(g(Y'))=r(Y')\leq n-1$. Therefore, $r(X \cup g(Y))\leq n$.

\section{ Distortion in groups of generalized piecewise-linear transformations 
}
This is the main section of this work. Here we introduce the groups we are going to deal with and we develop the necesary tools to prove \textbf{Theorem A}.

\begin{defn} \label{gp Rn}
\em{Let $n\in \mathbb{N}$ be any natural number. We define $\mathcal{R}_n$ as the set of all $f\in Homeo_+[0,1]$, such that there exist a sets 
(which we will call the set of break points and denote by $BP(f)$) such that $f$ is linear on each component of the complement of $BP(f)$ and 
$r(BP(f))\leq n$. Note that $\mathcal{R}_0=\PL_+ (I).$ }
\end{defn}

\begin{prop} \label{gp} 
\em{For each $n\in \mathbb{N}$, the set $\mathcal{R}_n$ is a group. }
\end{prop}

\noindent{\bf Proof.}  Let $f$ and $g$ be any elements of the set $\mathcal{R}_n$. Let $x$ be a point that belongs 
to both the complement of $g^{-1}(BP(f))$ and the complement of $BP(g)$. By the chain rule, we have 
\begin{align*}
    D(f\circ g)(x)=Df(g(x))\cdot Dg(x).
\end{align*}
We this conclude that $x$ belongs to the complement of $BP(f\circ g)$. This means that
$$BP(f\circ g) \subset g^{-1}(BP(f))\cup BP(g).$$
Since $r(BP(f))\leq n$ and $r(BP(g))\leq n$, Proposition \ref{prop rank} implies
\begin{align*}
    r(BP(f\circ g))\leq r(g^{-1}(BP(f))\cup BP(g)) \leq n. 
\end{align*}
This shows that $f\circ g$ belongs to $\mathcal{R}_n.$ Also, since $BP(f^{-1})=f^{-1}(BP(f))$, using Proposition \ref{prop rank} we have that $r(BP(f^{-1}))=r(BP(f))\leq n$. 
Therefore, $\mathcal{R}_n$ is a group.

\vspace{0.5cm}

We will refer to $\mathcal{R}_n$ as \textit{\textbf{group of generalized piecewise-linear transformations}}. 
By definition, we have that the groups $\mathcal{R}_n$ defined above give the following ascending chain of groups:
\begin{align*}
    \PL_+(I)=\mathcal{R}_0\subset\mathcal{R}_1\subset\mathcal{R}_2\subset\mathcal{R}_3\subset\dots
\end{align*}

In the next proposition, we build a length function in $\mathcal{R}_n$ that we will use to test that an element is undistorted. To do this, 
we will use the following notation: given a subset $A$ of the interval $[0,1]$, we denote by $A^{(n)}$ the subset that we obtain by applying 
derivation $n$ times to the set $A$ ({\em i.e.} $A^{(n)}$ is the $n^{\mathrm{th}}$ derived subset of $A$).

\begin{prop} \em{Given $n\in \mathbb{N}$, the function $L_n:\mathcal{R}_n\rightarrow \mathbb{N}$ defined as $L_n(f):=|BP(f)^{(n)}|$ 
is a length function. }
\end{prop}

\noindent{\bf Proof.} It is clear that $L_n(id)=0$ because the identity has no break point. Also, again by the chain rule, we have that
\begin{align*}
    BP(f)=f^{-1}(BP(f^{-1})).
\end{align*}
It thus follows from Proposition \ref{prop rank} that $L_n(f)=L_n(f^{-1})$.

It remains to prove the subadditive inequality. Let $f$ and $g$ be any elements in the group $\mathcal{R}_n$. We already know that
\begin{align*}
    BP(f\circ g) \subseteq g^{-1}(BP(f))\cup BP(g).
\end{align*}
Taking derived sets $n$ times we obtain finite sets, and 
\begin{align*}
   L_n(f\circ g) &=|BP(f\circ g)^{(n)}|\\
   &\leq |(g^{-1}(BP(f))\cup BP(g))^{(n)}|\\
    &=|g^{-1}(BP(f)^{(n)})\cup BP(g)^{(n)}|\\
    &=|g^{-1}(BP(f)^{(n)})|+|BP(g)^{(n)}|-|g^{-1}(BP(f)^{(n)})\cap BP(g)^{(n)}|\\
    &\leq   |g^{-1}(BP(f)^{(n)})|+|BP(g)^{(n)}|\\
    &=|(BP(f)^{(n)})|+| BP(g)^{(n)}|\\
    &=L_n(f)+L_n(g).
\end{align*}
We conclude that $L_n$ is a length function.

\begin{defn} 
\em{Let $n\in \mathbb{N}$ and let $f\in \mathcal{R}_n$ be any element. Since $(L_n(f^m))_{m\geq 0}$ is a subadditive sequence, by Fekete´s lemma, the following limit exists:
\begin{align*}
  l_n(f):=\lim_{m\rightarrow\infty} \frac{L_n(f^m)}{m}.  
\end{align*} }
\end{defn}

As before, this gives a tool to establish that an element in the group $\mathcal{R}_n$ is undistorted. The criterion is given by the proposition below.

\begin{prop} \label{no dist}
\em{Let $n\in \mathbb{N}$ and let $f\in \mathcal{R}_n$ be any element. If $l_n(f)>0$, then $f$ is undistorted in $\mathcal{R}_n$.}
\end{prop}

Let us recall the \textit{endpoint homomorphism}  $\eta:\PL_+[0,1]\rightarrow\mathbb{R}\oplus \mathbb{R}$ defined as $\eta (f)=(\log Df_+(0),\log Df_-(1)).$ To prove Theorem A, we will work with elements in the kernel of $\eta.$ 

\vspace{0.5cm}

\noindent{\textbf{Theorem A}}
{\em Let $n$ be a natural number. There exists $f\in Homeo_+[0,1]$ that is undistorted in $ \mathcal{R}_n$ but distorted in $\mathcal{R}_{n+1}$.}

\vspace{0.5cm}
\noindent{\bf Outline of the proof:} We will consider $f$ as a commutator in $\mathcal{R}_{n}\setminus\mathcal{R}_{n-1}$ that satisfies our criterion (Proposition \ref{no dist}). 
Then, we will deal with the distortion of $f$ in $\mathcal{R}_{n+1}$ with a combination of two techniques:
\begin{itemize}
    \item \textbf{(Diagonal trick)} The powers of the commutator $f$ can be written as the product of a bounded quantity of commutators. 
    \item \textbf{(Mather's argument)} This is a recipe to build functions (belonging to $\mathcal{R}_{n+1}$) allowing to write sequences of commutators with lengths in a prescribed sequence of natural numbers (see \cite{Avila}).
\end{itemize}

\noindent{\bf Building the element $f$.}
We fix $n\in \mathbb{N}$. Consider a sufficiently small interval $I$ that is contained $(0,1)$. Let $f_1$ satisfy the following properties:
\begin{itemize}
    \item We have that $supp(f_1)\subseteq I$ and $f_1(x)\geq x$ for all $x\in[0,1]$.
    \item There exist $x_0$ and $x_1$ in $Supp(f_1)$ such that $x_0<x_1$, $f_1(x_0)=x_1$ and $f_1$ is linear outside $[x_0,x_1]$.
    \item $f_1\in \mathcal{R}_n\setminus \mathcal{R}_{n-1}$.   
\end{itemize}

 Observe that $f_1$ cannot be distorted in $\mathcal{R}_n$. Indeed, for every $k\in \mathbb{N}$ we have $|BP(f_1^k)^{(n)}| =k\cdot|BP(f_1)^{(n)}|$, and since 
 $0<|BP(f_1)^{(n)}|<\infty$, this implies that $l_n(f_1)=|BP(f_1)^{(n)}|>0,$ which by Proposition \ref{no dist} shows that $f_1$ is undistorted in $\mathcal{R}_n.$\\ 

We want to build $f$ as a commutator. To do this, consider $t\in Ker(\eta)$ such that:
\begin{itemize}
    \item $t(x)\geq x$ for all $x\in [0,1].$
    \item $t(I)\cap I=\emptyset.$
\end{itemize}

Let $f_2:=tf_1t^{-1}$. We define $f:=[f_1,t]=f_1f_2^{-1}$. Note that $f_1$ and $f_2$ have disjoint support. Since $f_1$ is undistorted in $\mathcal{R}_n$, the element $f$ is also undistorted.

\vspace{0.3cm}

We next apply the diagonal trick to the element $f.$ This trick is well-known; we refers to \cite{Cal}, where in his study of the stable conmutator length in subgroups of $\PL_+(I)$, 
Calegari applies the diagonal trick in the group $\PL_+(I)$.

\vspace{0.2cm}

\noindent{\bf Diagonal trick:} We will prove that, for $m\geq 1$, the element $f^{m+1}$ can be written as a product of two commutators.

Let us take an interval $J$ contained in $(0,1)$ that contains the support of $f$. Let $h\in Ker(\eta)$ be such that:
\begin{itemize}
    \item $h(x)\geq x$ for all $x\in [0,1].$
    \item $h(J)\cap J=\emptyset.$ 
\end{itemize}

Observe that the last property implies that for every $i\in \mathbb{Z}\setminus\{0\}$, one has $h^i(J)\cap J=\emptyset.$ Consider the two finitely generated groups $G_0:=\langle f_1,t\rangle$ and $G_1:=\langle f_1,t,h\rangle$. The map $\Delta_m:G_0\rightarrow G_1$ defined as
 \begin{align*}
     \Delta_m(g):=\prod_{i=0}^{m}h^{i}gh^{-i},
 \end{align*}
 is an injective homomorphism. 

\begin{center}
\includegraphics[scale=0.5]{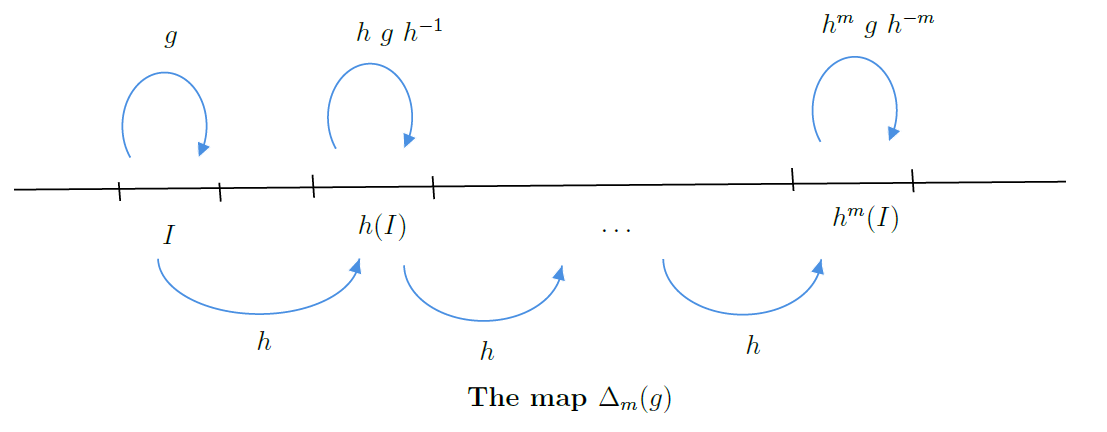}
\end{center}

 Let $f'\in G_1$ be defined as
 \begin{align*}
     f'=\prod_{i=0}^{m}h^{i}f^{i+1}h^{-i}.
 \end{align*}
We compute the commutator between $f'$ and $h$,
 \begin{align*}
     [f',h]&=f'\cdot (hf'^{-1}h^{-1})\\ 
     &=\left(\prod_{i=0}^{m}h^{i}f^{i+1}h^{-i}\right)\cdot \left(\prod_{i=1}^{m+1}h^{i}f^{-i}h^{-i}\right)\\
      &=\left(\prod_{i=0}^{m}h^{i}fh^{-i}\right)\cdot h^{m+1} f^{-(m+1)}h^{-(m+1)},
 \end{align*}
 where the last equality is justified because the maps in each product have support contained in one 
 of the sets $I, h(I),...,h^m(I)$ or $h^{m+1}(I)$, and maps with disjoint support commute. Then, we have the following relation:
 \begin{align*}
     [f',h]=\Delta_m(f) h^{m+1}f^{-(m+1)}h^{-(m+1)}.
 \end{align*}
 Finally, we can write $f^{m+1}$ as a product of two commutators:
 \begin{align*}
     f^{m+1}&=h^{-(m+1)}(\Delta_m(f)\cdot [f',h]^{-1}) h^{(m+1)}\\
     &=h^{-(m+1)}[\Delta_m(f_1),\Delta_m(t)]\,\,h^{(m+1)}\cdot h^{-(m+1)}[h,f']\,\, h^{(m+1)}\\
     &=[h^{-(m+1)}\Delta_m(f_1)\,h^{(m+1)},h^{-(m+1)}\Delta_m(t)\,h^{(m+1)}]\cdot [h,h^{-(m+1)}f'h^{(m+1)}],
 \end{align*}
where we use that the inverse of a commutator is a commutator, a homomorphism sends a commutator to a commutator,  
and conjugating a commutator gives us another commutator.

\vspace{0.5cm}

\noindent{\bf Mather's argument:} The next lemma is Mather's argument stated for our situation. 
This classical trick can be found (in different versions) 
in \cite{Avila}, \cite{CF}, \cite{Fis}, \cite{Le roux Kat}, \cite{Mil} and probably elsewhere. 
Let us point out that the statement is inspired by Avila's work (see \cite{Avila}).

\begin{lem} \label{avila trick}
{\em Let $[a,b]$ be an interval properly contained in (0,1) and let $n\in \mathbb{N}$. There exists a (linear) sequence of natural numbers $(k_m)$ such that for any sequences $(f_m)$ and $(g_m)$ of $\mathcal{R}_n\setminus \mathcal{R}_{n-1}$, such that $supp(f_m)\subseteq [a,b]$ and $supp(g_m)\subseteq [a,b]$, there exists a finite set $\mathcal{G}\subset \mathcal{R}_{n+1}$ such that in the subgroup of $\mathcal{R}_{n+1}$ generated by $\mathcal{G}$ one has $l_{\mathcal{G}}[f_m,g_m]\leq k_m.$ }
\end{lem}

\noindent{\bf Proof of the Lemma:} After conjugating by an element in $\PL_+[0,1] $, we can assume that $1/2$ belongs to $(a,b)$ and $0<b-a<1/2$. First we define three homeomorphisms in $\PL_+[0,1]$ that will be used along the proof. Let $0<a'<a<b<b'<1$. Let us define $h \in \PL_+[0,1]$ as $h^{-1}(x)=\frac{x}{2}+\frac{1}{4}$ on $[a',b']$ and being linear outside $[a',b'].$  Let $r\in$ $\PL_+[0,1]$ be defined as $r(x)=2x$ on $[0,\frac{3}{8}]$ and being linear on $[\frac{3}{8},1]$.

\begin{center}
\includegraphics[scale=0.5]{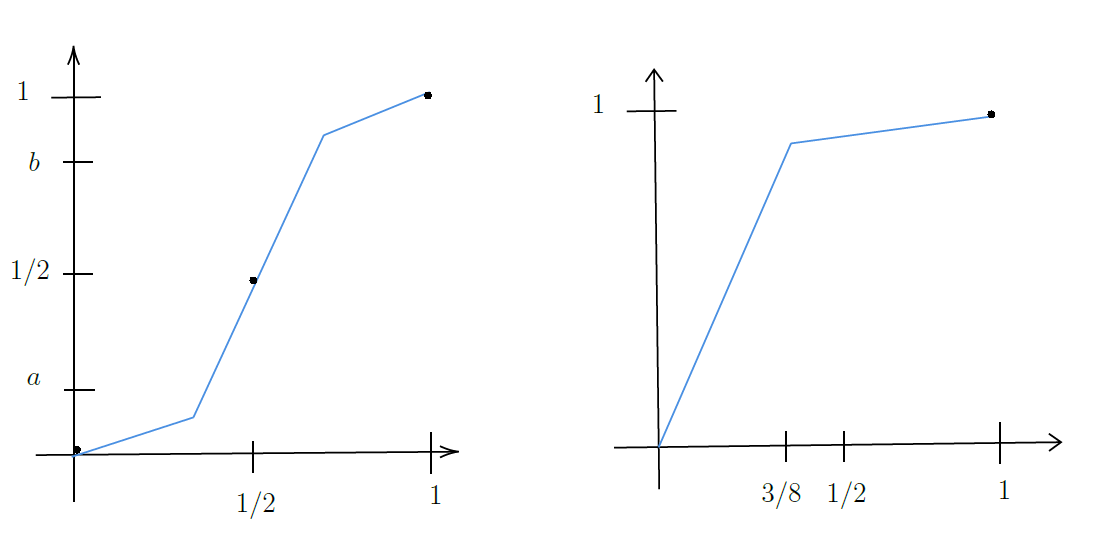}
\end{center}

Let $\alpha \in (0,1)$ be sufficiently small such that $\alpha+b<1$. We let $\widetilde{f} \in \ker(\eta)$ be such that
\begin{itemize}
    \item $\widetilde{f}(x)=x+\alpha$ for each $x\in [a,b]$.
    \item $supp(\widetilde{f})\subseteq [a',b']$ and $0<b'-a'<1/2$.
\end{itemize}

\begin{center}
\includegraphics[scale=0.5]{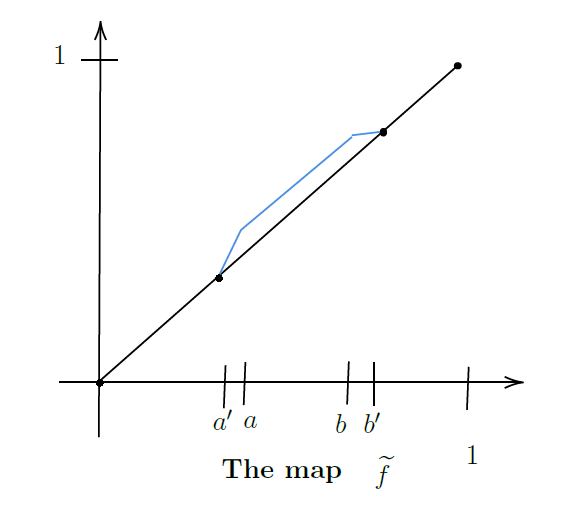}
\end{center}

We will define the sequence that we will use to compute commutators. Consider the sequence $(F_m)$ defined as
\begin{align*}
    F_m=r^{-m}h^{-m}\widetilde{f}h^mr^m.
\end{align*}
For each $m\geq1$, we denote $T_m$ the interval $ T_m:=[r^{-m}h^{-m}(a'),r^{-m}h^{-m}(b')]$. Note that $Supp(F_m)\subseteq T_m$. 

Let us prove that, for all $m\geq 1$, we have that $T_{m+1}\cap T_m=\emptyset.$ By construction, it is enough to show that $r^{-(m+1)}h^{-(m+1)}(b')<r^{-m}h^{-m}(a')$. An explicit computation shows that
\begin{align*}
T_m&=[r^{-m}h^{-m}(a'),r^{-m}h^{-m}(b')]\\
&=\left[\frac{a'}{2^{2m}}+\frac{1}{2^{2m+1}}+\frac{1}{2^{2m-1}}+\cdots+\frac{1}{2^{m+2}},\frac{b'}{2^{2m}}+\frac{1}{2^{2m+1}}+\frac{1}{2^{2m-1}}+\cdots+\frac{1}{2^{m+2}}\right].
\end{align*}
Since $0<b'-a'<1/2$, it follows that, for all $m\geq 1$, one has 
$$0<(2^{m+1}-3)-2(b'-a')+6a'.$$
This inequality implies that 
\begin{align*}
r^{-(m+1)}h^{-(m+1)}(b')&= \frac{b'}{2^{2m+2}}+\frac{1}{2^{2m+3}}+\frac{1}{2^{2m+2}}+\cdots+\frac{1}{2^{m+3}}\\ 
&<\frac{a'}{2^{2m}}+\frac{1}{2^{2m+1}}+\frac{1}{2^{2m-1}}+\cdots+\frac{1}{2^{m+2}}\\
&=r^{-m}h^{-m}(a'),
\end{align*}
and, therefore, $T_{m+1}$ and $T_m$ are disjoint sets.

\begin{center}
\includegraphics[scale=0.6]{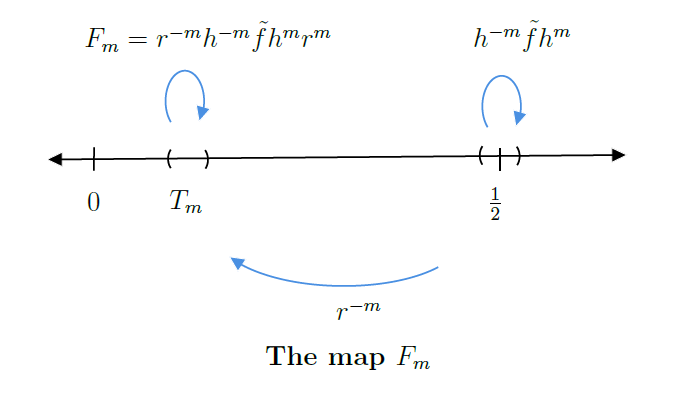}
\end{center}

We will build two elements $\widetilde{F}$ and $\widetilde{G}$ in $\mathcal{R}_{n+1}$ with support contained in $\cup_{m\geq1}T_m$. First, we define the set where the maps $\widetilde{F}$ and $\widetilde{G}$ will be supported. Let us fix $m_0\in \mathbb{N}$ such that 
\begin{equation} \label{m0}
    2^{m_0}\alpha >b-a.
\end{equation}
 For each $m\in \mathbb{N}$, we define the interval $T_m':=[r^{-m}h^{-m_0-m}(a),r^{-m}h^{-m_0-m}(b)]$. Let us justify the choice of $m_0$ because it implies that $T_m'$ and $F_m(T_m')$ are disjoint. Indeed, since for all $x\in [0,1]$ we have $F_m(x)\geq x$, to check this it is enough to prove that 
$$ r^{-m}h^{-m_0-m}(b)<F_m(r^{-m}h^{-m_0-m}(a)).$$ 
Using the definition of $F_m$, the previous inequality is equivalent to
$$h^{-m_0}(b)<\widetilde fh^m(a).$$
The last inequality is exactly the one in (\ref{m0}).  

 We define $\widetilde{F}$ as the map whose restriction to each $T_m'$ is given by
\begin{align*}
    \widetilde{F}|_{T_m'}=r^{-m}h^{-m_0-m}f_mh^{m_0+m}r^{m},
\end{align*} 
 and the identity elsewhere. Analogously, we build $\widetilde{G}$ using the sequence $(g_m)$. 

Observe that since the elements of each sequence $(f_m)$ and $(g_m)$ belong to $\mathcal{R}_n$, after taking the $n^{\mathrm{th}}$ 
derived subsets of $BP(\widetilde{F})$ and $BP(\widetilde{G})$, we obtain finite sets on each $T_m'.$ Since the intervals $T_m'$ accumulate at $0$, 
after taking the derived set once more, we obtain that $|BP(\widetilde{F})^{(n+1)}|=|BP(\widetilde{G})^{(n+1)}|=1$. Therefore, $\widetilde{F}$ and $\widetilde{G}$ belong to $\mathcal{R}_{n+1}$.

\begin{center}
\includegraphics[scale=0.5]{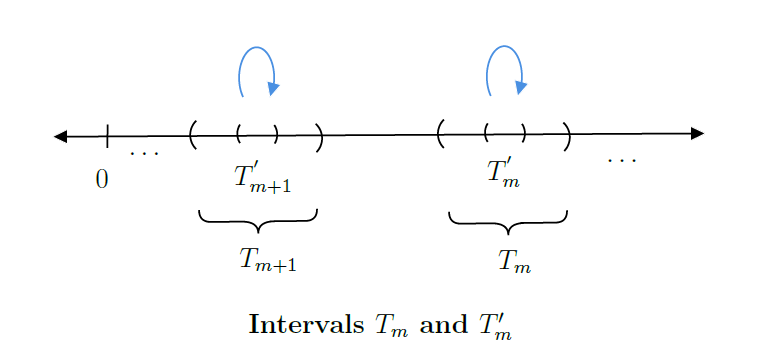}
\end{center}

We define $\mathcal{G}:=\{\widetilde{F},\widetilde{G},h,r,\widetilde{f}\} \subset \mathcal{R}_{n+1}$. For each $m\geq 1$, we consider the following elements in $\langle\mathcal{G}\rangle$: 
\begin{align*}
    A_m:=[\widetilde{F},F_m], \quad \,B_m:=[\widetilde{G},F_m],\quad \, C_m:=[\widetilde{F}^{-1} \widetilde{G}^{-1},F_m].
\end{align*}
We note that each of these homeomorphisms has support in the disjoint intervals $T_m'$ and $F_m(T_m')$. The maps $A_m,B_m$ and $C_m$ restricted to $T_m'$ are equal to $\widetilde{F}$, $\widetilde{G}$ and $\widetilde{F}^{-1}\widetilde{G}^{-1}$, respectively, and restricted to $F_m(T_m')$ are equal to $F_m\widetilde{F}^{-1}F_m^{-1} $, $F_m\widetilde{G}^{-1}F_m^{-1} $ and $F_m \widetilde{G}\widetilde{F} F_m^{-1} $, respectively. 

\begin{center}
\includegraphics[scale=0.5]{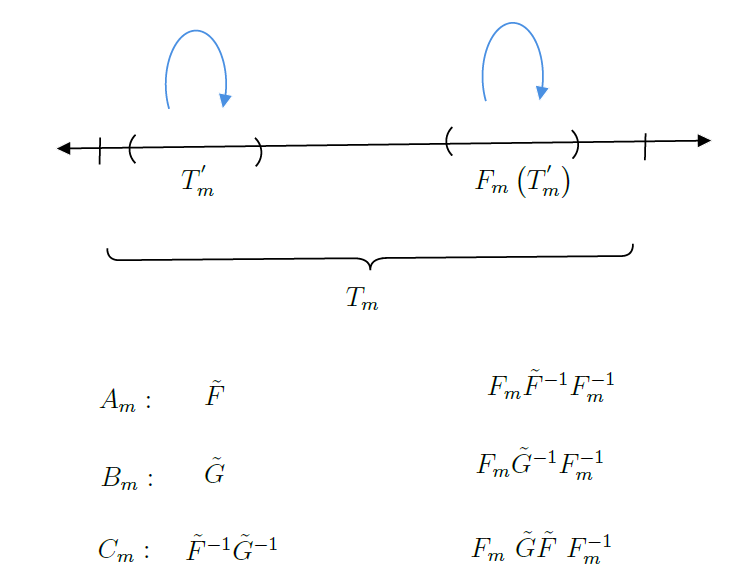}
\end{center}

Hence, $A_mB_mC_m$ has support in $T_m'$, and 
\begin{align*}
    A_mB_mC_m=r^{-m}h^{-m_0-m}f_mg_mf_m^{-1}g_m^{-1}h^{m_0+m}r^m.
\end{align*}
Then, $H_m:=h^{m_0+m}r^mA_mB_mC_m r^{-m}h^{-m_0-m} \in \langle \mathcal{G} \rangle$ satisfies that $H_m=[f_m,g_m]$ and $l_{\mathcal{G}}(H_m)\leq k_m=28m+2m_0+14$.

\vspace{0.32cm}

 \noindent{\textbf{Proof of Theorem A:} 
We consider the element $f$ built above. As we mentioned before, this element is undistorted in $\mathcal{R}_n.$

 Given the sequence $(k_m)$ of Lemma \ref{avila trick}, we fix a sequence $(i_m)$ of natural numbers such that $$\lim_{m\rightarrow\infty}\frac{k_m}{i_m+1}=0.$$
 By the Diagonal Trick, the powers of $f$ can be written as products of two commutators $f^{i_m+1}=[a_{i_m},b_{i_m}][c_{i_m},d_{i_m}]$ where 
 $$a_{i_m}:=h^{-(i_m+1)}\Delta_{i_m}(f_1)h^{(i_m+1)},\,\,\, b_{i_m}:=h^{-(i_m+1)}\Delta_{i_m}(t)\,h^{(i_m+1)},\,\,\,c_{i_m}:=h, $$
 $$d_{i_m}:=h^{-(i_m+1)}f'h^{(i_m+1)}.$$
 Also, from the construction we know that all the supports of the elements $a_{i_m},b_{i_m},c_{i_m}$ and $d_{i_m}$ are contained in the support of $h$, so there exists an interval properly contained in $(0,1)$ that contains the supports of $a_{i_m},b_{i_m},c_{i_m}$ and $d_{i_m}$. Then, applying Lemma \ref{avila trick} to the sequences of commutators, we obtain that there exist finite subsets $\mathcal{G}$ and $\mathcal{G'}$ of $\mathcal{R}_{n+1}$ such that $l_{\mathcal{G}}([a_{i_m},b_{i_m}])\leq k_m$ and $l_{\mathcal{G'}}([c_{i_m},d_{i_m}])\leq k_m.$  Hence, by the triangular inequality, 
\begin{align*}
    l_{\mathcal{G}\cup\mathcal{G'}}(f^{i_m+1})\leq 2 k_m,
\end{align*}
and therefore $f$ is distorted in $\mathcal{R}_{n+1}.$

\vspace{0.2 cm}

    {\em \noindent{\textbf{Question 4.}}  Is every element $f\in$ $\PL_+(I)$  distorted in $\mathcal{R}_1$? More generally, is every element of $\mathcal{R}_n$ distorted in $\mathcal{R}_{n+1}$? }
\vspace{0.2 cm}

\begin{obs}
{Let us define for each $n\in \mathbb{N}$ the following set     $$\widetilde{\mathcal{R}_n}:=\{f\in \mathcal{R}_n : f \mbox{ is bilipchitz}\}.$$
    With this definition, again we obtain a chain of groups as before
    \begin{align*}
\PL_+(I)\leq\widetilde{\mathcal{R}_1}\leq\widetilde{\mathcal{R}_2}\leq\widetilde{\mathcal{R}_3}\leq\dots
\end{align*}
Similarly, we could ask whether there exist a distorted elements in $\widetilde{\mathcal{R}_{n+1}}$ that is undistorted in $\widetilde{\mathcal{R}_{n}}$. 
However, the situation here is different, because  there is no distorted elements in the groups $\widetilde{\mathcal{R}_{n}}$. Indeed, we can consider 
the function $L:\widetilde{\mathcal{R}_{n}}\rightarrow\mathbb{R}$ defined as
    \begin{align*}
        L(f):=\log (||f||_{Lip}),
    \end{align*}
    where $||\cdot||_{Lip}$ is taken in the complement of the set of break points. We see that $L$ is a length function and if $f\in \widetilde{\mathcal{R}_{n}}$ is different of the identity, then we have that $L(f)>0$ and $L(f^n)\geq n\cdot L(f).$ Therefore, $f$ is undistorted.}
\end{obs}

\begin{rem}
    Let us explain why the combination of the Diagonal Trick and the Mather's argument does not work for the groups $\widetilde{\mathcal{R}_{n}}$. Indeed, as $m$ increases, the commutators of the diagonal arguments take account arbitrarily large powers of the element $f$, and so the elements ($\widetilde{F}$ and $\widetilde{G}$) created by Mather's trick cannot belong to $\widetilde{\mathcal{R}_{n+1}}$ because they are not bilipchitz. 
\end{rem}

\vspace{0.35cm}

\noindent{\bf Acknowledgments.} I would like to gratefully acknowledge to Yash Lodha for suggesting the problem  
and Andrés Navas for proposing it to me as well as for the guidance and many encouragements during the work. 
Also, I thank to Hélène Eynard-Bontemps, Cristóbal Rivas and Michelle Triestino for the comments on previous 
versions of this article. I thank the hospitality 
of IHP during the trimester "Around the Zimmer program", I acknowledge support of the Institut Henri Poincaré 
(UAR 839 CNRS-Sorbonne Université) and LabEx CARMIN (ANR-10-LABX-59-01). I was partially funded by Fondecyt regular 1220032.

%%%%%%%%%%%%%%%%%%%%%%%%%%%%%%%%%%%%%%%%%%%%%%%%%%%%%%%%%%%%%%%%%%%%%%%%

\begin{small}

\vspace{0.3cm}

\noindent Leonardo Dinamarca Opazo

\noindent Dpto. de Matem\'atica y Ciencias de la Computación, Univ. de Santiago de Chile (USACH)

\noindent Alameda Lib. Bdo O'Higgins 3363, Estaci\'on Central, Santiago, Chile

\noindent Emails: leonardo.dinamarca@usach.cl

\end{small}

\end{document}